\documentclass[sn-basic]{sn-jnl}

\usepackage{graphicx}%
\usepackage{multirow}%
\usepackage{amsmath,amssymb,amsfonts}%
\usepackage{amsthm}%
\usepackage{mathrsfs}%
\usepackage[title]{appendix}%
\usepackage{xcolor}%
\usepackage{textcomp}%
\usepackage{manyfoot}%
\usepackage{booktabs}%
\usepackage{algorithm}%
\usepackage{algorithmicx}%
\usepackage{algpseudocode}%
\usepackage{listings}%
\usepackage{graphicx}
\usepackage{epstopdf}
\begin{document}

\title[Article Title]{Distinguishing between long-transient and asymptotic states in a biological aggregation model}


\author*[1]{\fnm{Jonathan R.} \sur{Potts}}\email{j.potts@sheffield.ac.uk}

\author[2]{\fnm{Kevin J.} \sur{Painter}}\email{kevin.painter@polito.it}
			
\affil*[1]{\orgdiv{School of Mathematics and Statistics}, \orgname{University of Sheffield}, \orgaddress{\street{Hounsfield Road}, \city{Sheffield}, \postcode{S3 7RH}, \country{UK}}}

\affil[2]{\orgdiv{Dipartimento Interateneo di Scienze}, \orgname{Progetto e Politiche del Territorio (DIST)}, \orgaddress{\city{ Politecnico di Torino}, \country{Italy}}}

\abstract{Aggregations are emergent features common to many biological systems. Mathematical models to understand their emergence are consequently widespread, with the aggregation-diffusion equation being a prime example.  Here we study the aggregation-diffusion equation with linear diffusion.  This equation is known to support solutions that involve both single and multiple aggregations.  However, numerical evidence suggests that the latter, which we term `multi-peaked solutions' may often be long-transient solutions rather than asymptotic steady states.  We develop a novel technique for distinguishing between long transients and asymptotic steady states via an energy minimisation approach.  The technique involves first approximating our study equation using a limiting process and a moment closure procedure.  We then analyse local minimum energy states of this approximate system, hypothesising that these will correspond to asymptotic patterns in the aggregation-diffusion equation.  Finally, we verify our hypotheses through numerical investigation, showing that our approximate analytic technique gives good predictions as to whether a state is asymptotic or a long transient.  Overall, we find that almost all twin-peaked, and by extension multi-peaked, solutions are transient, except for some very special cases.  We demonstrate numerically that these transients can be arbitrarily long-lived, depending on the parameters of the system.}




\maketitle

	\section{Introduction}

Aggregation phenomena are widespread in biology, from cell aggregations \citep{budrene1995} to the swarming \citep{roussi2020}, schooling \citep{makris2009}, flocking \citep{clark1984foraging}, and herding \citep{bond2019fission} of animals.  When modelled from a continuum perspective (as opposed to via interacting particles), the principal tools take the form of partial differential equations with non-local advection, sometimes combined with a diffusive term \citep{topazetal2006}.  Indeed, such equations are often called aggregation equations \citep{laurent2007local}, highlighting their importance in modelling aggregations, or aggregation-diffusion equations \citep{carrillo2019aggregation} if there is a diffusion term.  

As well as modelling aggregated groups of organisms, such equations have also been used to model aggregation-like phenomena elsewhere, such as animal home ranges and territories \citep{briscoeetal2002, pottslewis2016a} and consensus convergence in opinion dynamics \citep{garnier2017consensus}.  This very broad range of applications, together with the mathematical complexity in dealing with nonlinear nonlocal partial differential equations (PDEs), has led to a great amount of interest from applied mathematicians in understanding the properties of these PDEs \citep{painter2023biological}.  

Of particular interest from a biological perspective are the pattern formation properties of aggregation-diffusion equations, since these can reveal the necessary processes required for observed patterns to emerge.  Many traditional techniques for analysing pattern formation, such as linear stability analysis and weakly nonlinear analysis, focus on the onset of patterns from small perturbations of a non-patterned (i.e. spatially homogeneous) state.  However, patterns observed in actual biological systems will often be far from the non-patterned state, and not necessarily emerge from small perturbations of spatially homogeneous configurations \citep{krause2020one, veerman2021beyond}.  

Sometimes observed patterns will be asymptotic steady states or other types of attractors.  But frequently biological systems will be observed in transient states \citep{hastings2018transient,morozov2020long}.  These transient states may persist for a very long time, sometimes so long that they are hard to distinguish from asymptotic states.  Moreover, as well as transients being difficult to decipher from observations of biological systems, they can also be tricky to determine from numerical solutions of a PDE model.  Therefore analytic techniques are required to guide those engaging in numerical analysis of PDEs as to whether the solution they are observing is likely to be a long transient or an asymptotic state.

Our aim here is to provide such analytic techniques for a class of 1D aggregation-diffusion equations of the following form  
\begin{align}
	\label{eq:agg}
	\frac{\partial u}{\partial t}=D\frac{\partial^2 u}{\partial x^2}-\gamma\frac{\partial}{\partial x}\left[u\frac{\partial}{\partial x}(K\ast u)\right],
\end{align}
where $K$ is a non-negative averaging kernel, symmetric about $0$, with $\|K\|_{L^\infty}=1$, and 
\begin{align}
	\label{eq:conv}
	K\ast u(x)=\int_{\Omega}K(z)u(x+z){\rm d}z
\end{align}
is a convolution, where $\Omega$ is the spatial domain of definition.  Here, $D$ and $\gamma$ are constants, and $\Omega$ is the circle given by interval $[-L,L]$ with periodic boundary conditions imposed.

Our approach is not exact, in the sense that we approximate our study PDE first through the limit as $D/\gamma \rightarrow 0$, then via a moment closure assumption.  However, this approximation allows us to analyse the associated energy functional, finding explicit mathematical expressions for local energy minima.  Our conjecture is that local energy minima of the approximate system are qualitatively similar to the asymptotic patterns observed the aggregation-diffusion equation we are studying, but any states that do not represent local energy minima of the approximate system are transient states.  We then test this numerically in some specific cases.  

Of particular interest is the question of whether multi-peaked solutions are asymptotic steady states or long transients, which is the question that originally motivated this work.  Various numerical studies of Equation (\ref{eq:agg}), and similar equations, report multi-peaked solutions \citep{armstrong2006continuum, buttenschon2020non, carrillo2019aggregation, daneri2022deterministic}.  However, merging and decaying of peaks have also been observed.  Furthermore, analytic investigations into chemotaxis equations, which have some similarities with aggregation equations, have demonstrated that multi-peaked solutions can often be long transients \citep{potapov2005metastability}.  

This work demonstrates that, except for the very specific case where peaks are of identical heights and evenly-spaced, any two-peaked solutions will eventually evolve into a solution with at most one peak, as the smaller peak decays to zero. The time it takes for the smaller peak to decay grows rapidly with the start height of the smaller peak, eventually tending to infinity as the difference in start heights between the two peaks tends to zero.  We show that a key parameter governing the speed of this decay is the diffusion constant $D$, with higher diffusion constants leading to faster decays.  We conjecture that, as $D \rightarrow 0$, the time to decay tends to infinity, meaning that two-peaked solutions become  stable.

Finally, we investigate the effect of incorporating logistic growth of the population into our model.  The motivation for this is that, in situations where transient solutions exists for a long time, it is no longer biologically reasonable to assume that we are working in situations where births and deaths are negligible.  We show that, for a given set of parameters and initial condition, there is a critical net reproduction rate, below which the smaller peak will decay and above which it will persist.

\section{Methodological approach}
\label{sec:mod}

\begin{figure}[t]
	\begin{center}
		\includegraphics[width=\columnwidth]{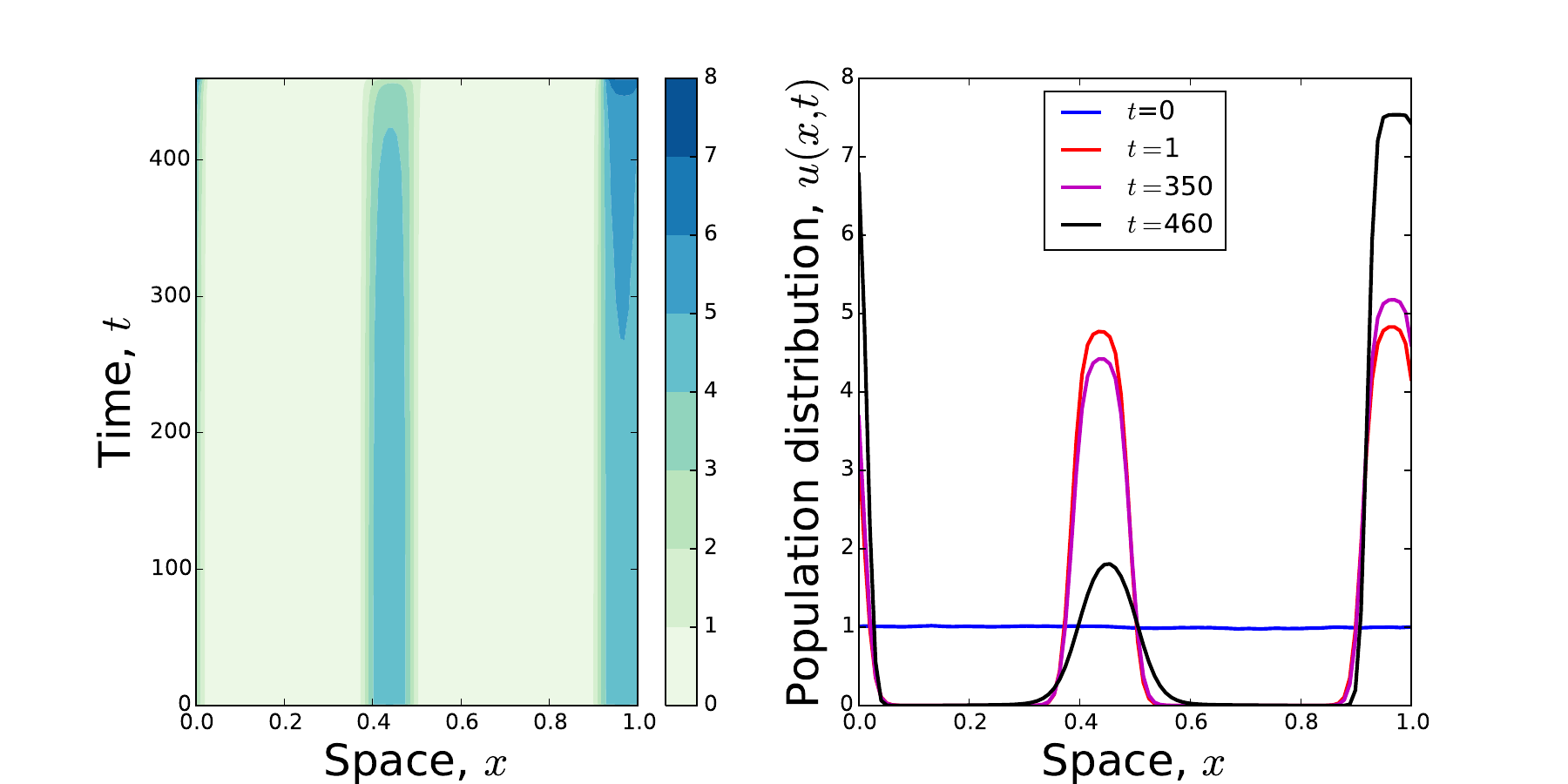}
	\end{center}
	\caption{Numerical solutions of Equation (\ref{eq:agg}) starting with initial conditions that are a small random fluctuation of the constant steady state.  By $t=1$ clear aggregations have formed that might seem stable were the solution only run to around time $t=100$.  However, if we run the solution further in time, we see that the middle peak is gradually decaying, and this decay is speeding up over time, so that by $t=460$ the peak in the middle is much smaller than the other peak.  Here, $D=1$, $\gamma=10$, and $K$ is a top-hat kernel (Equation \ref{eq:Kth}) with $\delta=0.1$.}
	\label{fig:nle_rand_ic}
\end{figure}

Our study is motivated by an observation.  Often, when simulating Equation (\ref{eq:agg}), multiple aggregations may form and persist for a very long time.  This can give the appearance of multi-peaked asymptotic stable states.  For example, Figure \ref{fig:nle_rand_ic} shows a numerical solution where two peaks have formed by time $t=1$.  These appear stable on timescales up to two orders of magnitude longer than the time they took to form: even by time $t=100$, the solution has not changed very much (Figure \ref{fig:nle_rand_ic}a).  However, if we keep running the simulation, we see one of the peaks decay and the other slowly swallow up the former's mass.  The question then arises whether multi-peaked solutions to Equation (\ref{eq:agg}) are ever actually stable, or are they always just long transients?

To answer this question, our approach will not be to analyse Equation (\ref{eq:agg}) directly, but rather take two approximations, which enable us to perform analytic calculations.  First, we assume that $\gamma \gg D$.  Second, we make the following moment closure assumption
\begin{align}
	\label{eq:conv2}
	K\ast u(x) &\approx u+\frac{\sigma^2}{2}\frac{\partial^2 u}{\partial x^2}
\end{align}
where
\begin{align}
	\label{eq:conv3}
	\sigma^2=\int_{-L}^L x^2K(x){\rm d}x
\end{align}
is the second moment of $K$.  This leads to the following approximate version of Equation (\ref{eq:agg})
\begin{align}
	\label{eq:agg2}
	\frac{\partial u}{\partial t}=-\gamma\frac{\partial}{\partial x}\left[u\left(\frac{\partial u}{\partial x}+\frac{\sigma^2}{2}\frac{\partial^3 u}{\partial x^3}\right)\right].
\end{align}
Note that Equations (\ref{eq:agg}) and (\ref{eq:agg2}) both preserve mass when solved with periodic boundary conditions (i.e. $u(x,L)=u(x,-L)$ and $\frac{\partial u}{\partial x}(x,L)=\frac{\partial u}{\partial x}(x,-L)$), so that if we define
\begin{align}
	\label{eq:p0}
	p:=\int_{-L}^{L}u(x,0){\rm d}x
\end{align}
then
\begin{align}
	\label{eq:pt}
	\int_{-L}^{L}u(x,t){\rm d}x = p,
\end{align}
for all $t>0$.

Our tactic will be to search for minimum energy solutions to Equation (\ref{eq:agg2}) using the following energy functional
\begin{align}
	\label{eq:energy}
	E[u]=-\int_{-L}^L u\left(u+\frac{\sigma^2}{2}\frac{\partial^2 u}{\partial x^2}\right){\rm d}x.
\end{align}
In particular, we are interested in examining critical points of $E[u]$, so calculate
\begin{align}
	\frac{\partial E}{\partial t}&=-\int_{-L}^L \left[\frac{\partial u}{\partial t}\left(u+\frac{\sigma^2}{2}\frac{\partial^2 u}{\partial x^2}\right)+u\left(\frac{\partial u}{\partial t}+\frac{\sigma^2}{2}\frac{\partial^2 }{\partial x^2}\frac{\partial u}{\partial t}\right)\right] {\rm d}x \nonumber \\
	&=-\int_{-L}^L 2\frac{\partial u}{\partial t}\left(u+\frac{\sigma^2}{2}\frac{\partial^2 u}{\partial x^2}\right) {\rm d}x \nonumber \\
	&=2\gamma\int_{-L}^L \frac{\partial}{\partial x}\left[u\frac{\partial }{\partial x}\left(u+\frac{\sigma^2}{2}\frac{\partial^2 u}{\partial x^2}\right)\right]\left(u+\frac{\sigma^2}{2}\frac{\partial^2 u}{\partial x^2}\right) {\rm d}x \nonumber \\
	&=-2\gamma\int_{-L}^L u\frac{\partial }{\partial x}\left(u+\frac{\sigma^2}{2}\frac{\partial^2 u}{\partial x^2}\right)\frac{\partial }{\partial x}\left(u+\frac{\sigma^2}{2}\frac{\partial^2 u}{\partial x^2}\right) {\rm d}x \nonumber \\
	&=-2\gamma\int_{-L}^L u\left[\frac{\partial }{\partial x}\left(u+\frac{\sigma^2}{2}\frac{\partial^2 u}{\partial x^2}\right)\right]^2 {\rm d}x.
	\label{eq:energy_dec}
\end{align}
Here, the second and fourth equalities use integration by parts, together with the periodic boundary conditions.  
If we assume that there exist non-negative solutions to Equation (\ref{eq:agg2}) then the final expression in Equation (\ref{eq:energy_dec}) is non-positive, so that $E[u]$ is non-increasing. Whilst we do not currently have a proof of the non-negativity of $u$, we note that that all our numerics suggest that non-negativity is preserved over time, that non-negativity results exist for Equation (\ref{eq:agg}) for a variety of different kernels $K$ \citep{carrillo2019aggregation, giunta2022local, jungel2022nonlocal}, and so conjecture these might be transferable to the situation of Equation (\ref{eq:agg2}) with some effort.

Equation (\ref{eq:energy_dec}) shows that critical points, $u_*(x)$, of the energy functional occur when
\begin{align}
	\int_{-L}^L u_*\left[\frac{\partial }{\partial x}\left(u_*+\frac{\sigma^2}{2}\frac{\partial^2 u_*}{\partial x^2}\right)\right]^2 {\rm d}x = 0,
\end{align}
which means that, on any connected subset of $[-L,L]$, either $u_*(x)=0$ or 
\begin{align}
	&u_*+\frac{\sigma^2}{2}\frac{\partial^2 u_*}{\partial x^2} = C \nonumber \\
	\implies&u_*(x)=C+A\sin\left(\frac{x\sqrt{2}}{\sigma}\right)+B\cos\left(\frac{x\sqrt{2}}{\sigma}\right)
	\label{eq:crit_energy}
\end{align}
for constants $A$, $B$, and $C$.

Numerics suggest that Equation (\ref{eq:agg}) tends towards a solution containing one or many aggregations, interspersed by constant sections close or near to zero (e.g. Figure \ref{fig:nle_rand_ic}).  We want to construct differentiable solutions that have this type of qualitative appearance, yet also correspond to critical points of $E[u]$.  These can be constructed piecewise from Equation (\ref{eq:crit_energy}).  For example, as long as $\pi\sigma<\sqrt{2}L$, a single-peaked solution can be given as follows
\begin{align}
	u_*(x)=\begin{cases}
		\epsilon+c_\epsilon\left[1+\cos\left(\frac{x\sqrt{2}}{\sigma}\right)\right], & \mbox{if }x \in \left(-\frac{\pi\sigma}{\sqrt{2}},\frac{\pi\sigma}{\sqrt{2}}\right) \\
		\epsilon,  & \mbox{otherwise,}
	\end{cases}
	\label{eq:one_hump}
\end{align}
where $\epsilon \in \left[0,\frac{p}{2L}\right]$ and $c_\epsilon$ are constants.  One can also construct multi-peaked solutions in a similar way (which we will do later in the case of two peaks).  Notice that such solutions are continuously differentiable, i.e. $u_* \in C^1([-L,L])$, but not necessarily twice differentiable, so need to be understood in a weak sense \citep{evans2022partial}.  

By Equation (\ref{eq:pt}), a direct calculation gives
\begin{align}
	\label{eq:cep}
	c_\epsilon=\frac{p-2\epsilon L}{\sqrt{2}\pi\sigma}
\end{align}
so that the only free parameter in Equation (\ref{eq:one_hump}) is $\epsilon$.  Since the energy, $E[u]$, is non-increasing over time, the question arises as to which value of $\epsilon$ minimises $E[u]$ across the set of all functions of the form in Equation (\ref{eq:one_hump}).  Our approach is to derive such minima, both in the example from Equation (\ref{eq:one_hump}) and in various multi-peaked examples, conjecturing that such minima ought to approximate asymptotic solutions to the original problem in Equation (\ref{eq:agg}).  We then test these conjectures by investigating Equation (\ref{eq:agg}) numerically.  
\begin{figure}[t]
	\begin{center}
		\includegraphics[width=\columnwidth]{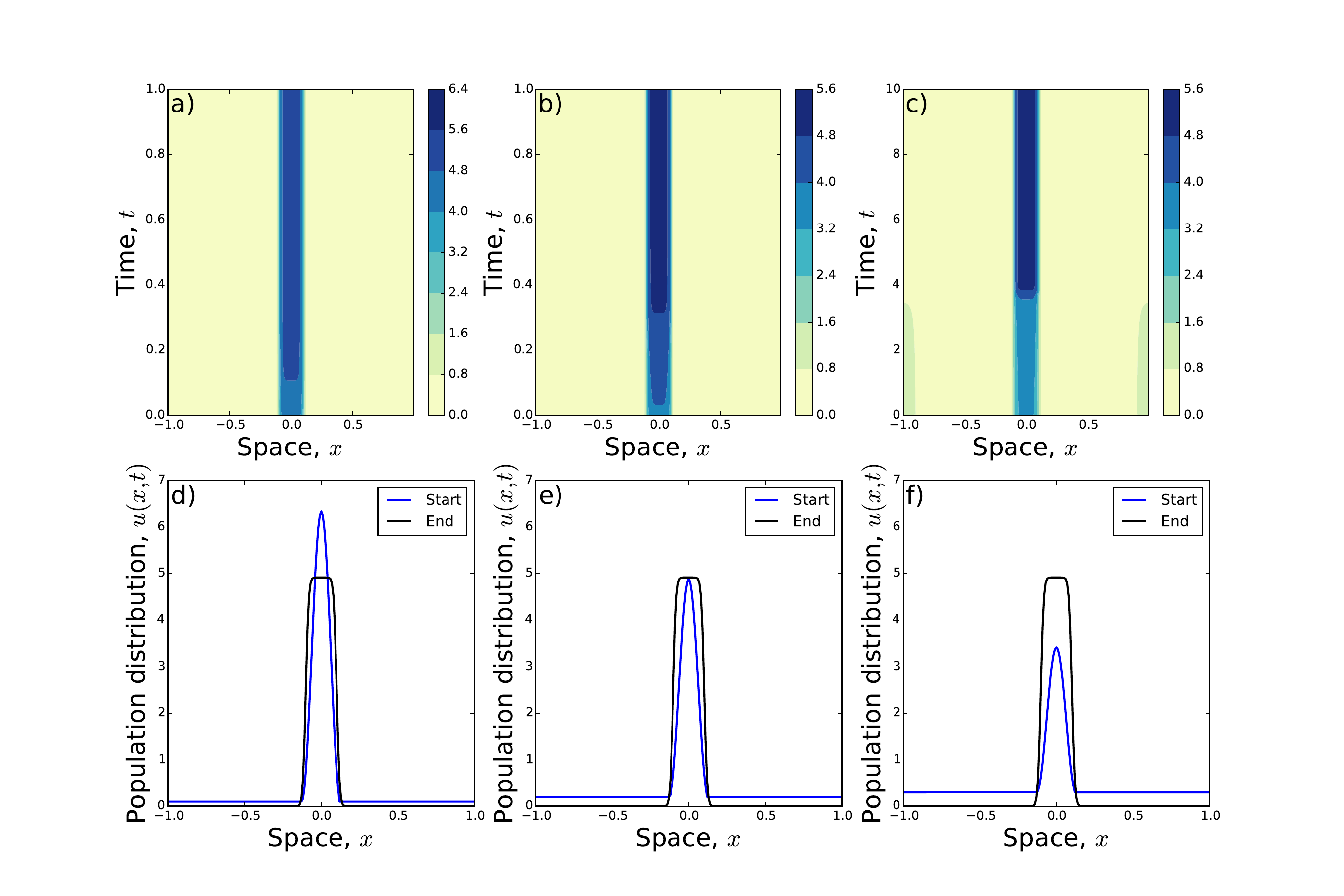}
	\end{center}
	\caption{When the initial condition is a single peak surrounded by an area of constant density $\epsilon$, that area becomes sucked-up into the peak. Panels (a) and (d) show this for $\epsilon=0.1$; (b) and (e) have $\epsilon=0.2$; (c) and (f) have $\epsilon=0.3$.  In the latter case, a second peak emerges at $x=\pm 1$ but decays by around $t \approx 4$, to leave a single-peaked final state.  Panels (a-c) show the time-evolution of the system.  Panels (d-f) show the initial conditions (blue curves) and final states (black). In all panels, $D=1$, $\gamma=10$, and $K$ is a top-hat kernel (Equation \ref{eq:Kth}) with $\delta=0.1$.}
	\label{fig:nle_single_peak}
\end{figure}

\section{Single peak}
\label{sec:one_peak}


Combining Equations (\ref{eq:energy}) and (\ref{eq:one_hump}) gives
\begin{align}
	\label{eq:energy1p}
	E[u_*]=-\int_{-L}^L u_*\left(u_*+\frac{\sigma^2}{2}\frac{{\rm d}^2 u_*}{{\rm d} x^2}\right){\rm d}x.
\end{align}
Now, for $-\pi\sigma/\sqrt{2}<x<\pi\sigma/\sqrt{2}$, we have that 
\begin{align}
	u_*(x)=\epsilon+c_\epsilon\left[1+\cos\left(\frac{x\sqrt{2}}{\sigma}\right)\right]
	\label{eq:one_hump_a}
\end{align}
which is a solution to 
\begin{align}
	u_*+\frac{\sigma^2}{2}\frac{\partial^2 u_*}{\partial x^2} = \epsilon+c_\epsilon.
\end{align}
Hence 
\begin{align}
	\label{eq:energy1p2}
	E[u_*]&=-\int_{-\frac{\pi\sigma}{\sqrt{2}}}^{\frac{\pi\sigma}{\sqrt{2}}}\left[\epsilon+c_\epsilon\left(1+\cos\left(\frac{\sqrt{2}x}{\sigma}\right)\right)\right](\epsilon+c_\epsilon){\rm d}x-2\int_{\frac{\pi\sigma}{\sqrt{2}}}^L\epsilon^2{\rm d}x
	&=-\pi\sigma\sqrt{2}(c_\epsilon^2+2\epsilon c_\epsilon)-2L\epsilon^2.
\end{align}
Using Equation (\ref{eq:cep}) and rearranging gives
\begin{align}
	\label{eq:energy1p3}
	E[u_*]=\frac{2L}{\pi\sigma}(\pi\sigma-\sqrt{2}L)\epsilon^2 + \frac{2p}{\pi\sigma}(\sqrt{2}L-\pi\sigma)\epsilon - \frac{p^2}{\sqrt{2}\pi\sigma}.
\end{align}
Since $\pi\sigma<\sqrt{2}L$ (see above Equation \ref{eq:one_hump}), this is a negative quadratic in $\epsilon$.  Furthermore, the maximum is where $\epsilon=\frac{p}{2L}$.  Now, $\epsilon \in \left[0,\frac{p}{2L}\right]$, so $E[u_*]$ is an increasing function of $\epsilon$ on the interval $\left[0,\frac{p}{2L}\right]$.  Hence the minimum energy is where $\epsilon=0$.

This analysis suggests that if a numerical solution to either Equation (\ref{eq:agg}) or (\ref{eq:agg2}) results in a single peak at long times, we might expect that peak to be of a similar form to Equation (\ref{eq:one_hump}) with $\epsilon=0$.  We test this conjecture by solving Equation (\ref{eq:agg}) numerically with initial conditions given by Equation (\ref{eq:one_hump}) for various different values of $\epsilon \in \left[0,\frac{p}{2L}\right]$, fixing $p=L=1$.  For these simulations, we set $D=1$, $\gamma=10$, and
\begin{align}
	\label{eq:Kth}
	{K}(x)=\begin{cases}\frac{1}{2\delta} & \mbox{for $-\delta<x<\delta$} \\
		0 & \mbox{otherwise,}
	\end{cases}
\end{align}
so that $\sigma=\delta/\sqrt{3}$.  Numerics reveal that the system does indeed tend towards a single-peaked solution, where the width of the peak is approximately $\sqrt{2}\pi\sigma$ and the solution is zero elsewhere (Figure \ref{fig:nle_single_peak}).  However, the asymptotic distribution is more flat-topped than the initial condition, owing to the fact that the initial condition arises from a moment closure approximation of $K\ast u$.  This approximation reduces the analytic solution to a single Fourier mode, whereas the numerical solution could have arbitrarily many Fourier modes. 

Finally note that, in the case $\epsilon=0.3$ (Figure \ref{fig:nle_single_peak}c,f), a second peak emerges around $x=\pm 1$ (which are identified due to the periodic boundaries, recalling that $L=1$).  However, this decays by about $t=4$.  We will return to this phenomenon of decaying secondary peaks in the next section.

\section{Twin peaks}
\label{sec:twin_peaks}

In this section, we examine situations where there are two peaks.  First, we look at situations where the peaks are  the same height, then at cases where one peak is smaller than the other.  

\subsection{Peaks of identical height}
\label{sec:same_height}

Similar to the single-peak case, here we want to understand whether it is energetically favourable for a solution to have no mass outside the two peaks.  More precisely, we examine the energy of the following solution to Equation (\ref{eq:agg2}), which is a critical point of $E[u]$
\begin{align}
	u_*(x)=\begin{cases}
		\epsilon+c_{\epsilon}\left[1+\cos\left(\frac{(x+x_0)\sqrt{2}}{\sigma}\right)\right], & \mbox{if }x \in \left(-x_0-\frac{\pi\sigma}{\sqrt{2}},-x_0+\frac{\pi\sigma}{\sqrt{2}}\right), \\
		\epsilon+c_{\epsilon}\left[1+\cos\left(\frac{(x-x_0)\sqrt{2}}{\sigma}\right)\right], & \mbox{if }x \in \left(x_0-\frac{\pi\sigma}{\sqrt{2}},x_0+\frac{\pi\sigma}{\sqrt{2}}\right), \\
		\epsilon,  & \mbox{otherwise.}
	\end{cases}
	\label{eq:two_humps1}
\end{align}
Here, $x_0 \in \left(\frac{\pi\sigma}{\sqrt{2}},\frac{L}{2}\right)$ is half the (shortest) distance between the centres of the two peaks. As in the single-peak case, we can use Equation (\ref{eq:pt}) to calculate 
\begin{figure}[t]
	\begin{center}
		\includegraphics[width=\columnwidth]{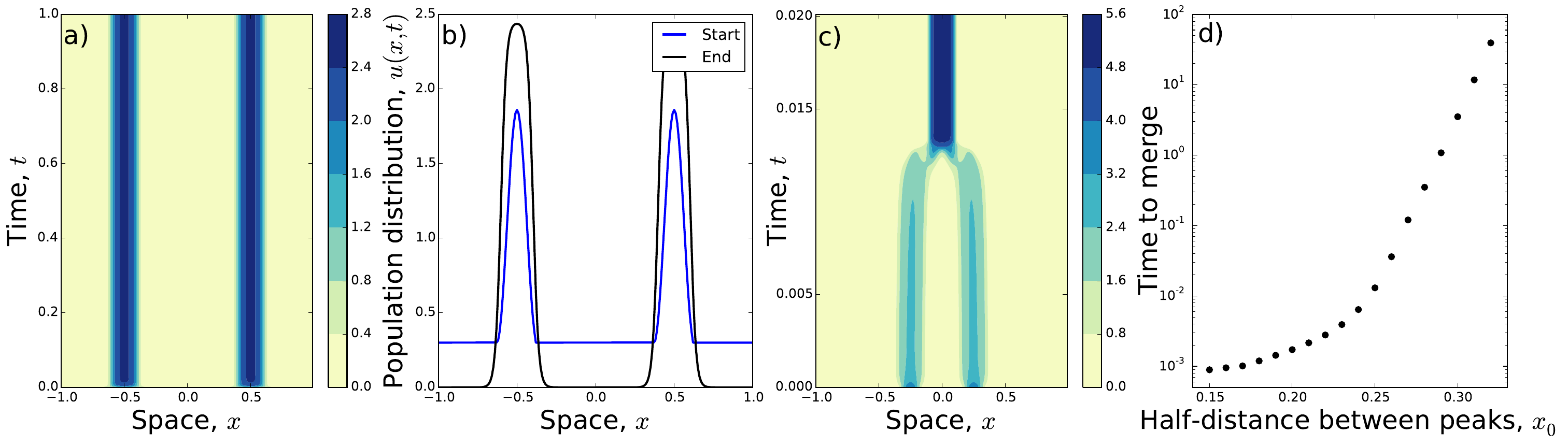}
	\end{center}
	\caption{Similar to the single peak case (Figure \ref{fig:nle_single_peak}), when we start with two peaks of equal heights, surrounded by an area of constant density $\epsilon$, that area becomes sucked-up into the peak. Panels (a) and (b) show this for $x_0=0.5$ and $\epsilon=0.2$, where both peaks remain.  For $x<0.5$, peaks merge, shown in Panel (c) for $x_0=0.25$. Panel (d) shows the time to merge as a function of $x_0$. Parameters $D$, $\gamma$, and $K$ are as in Figure \ref{fig:nle_single_peak}.}
	\label{fig:nle_twin_peaks}
\end{figure}
\begin{align}
	c_{\epsilon}=\frac{p-2L\epsilon}{2\sqrt{2}\pi\sigma}.
\end{align}
A direct calculation using the definition of $E[u]$ from Equation (\ref{eq:energy}) leads to
\begin{align}
	\label{eq:energy2a}
	E[u_*]=\frac{\sqrt{2}L}{\pi\sigma}(\sqrt{2}\pi\sigma-L)\epsilon^2 + \frac{\sqrt{2}p}{\pi\sigma}(L-\sqrt{2}\pi\sigma)\epsilon - \frac{p^2}{2\sqrt{2}\pi\sigma}.
\end{align}
Since $\sqrt{2}\pi\sigma<L$, this is a negative quadratic in $\epsilon$.  The unique turning point is a maximum at $\epsilon=\frac{p}{2L}$, so $E[u_*]$ is an increasing function of $\epsilon$ on the interval $\left[0,\frac{p}{2L}\right]$.  Hence the minimum energy in the two-peaked case is where $\epsilon=0$, as with the one-peaked case.  However, comparing the $\epsilon=0$ situation with one peak (Equation \ref{eq:energy1p3}), against that with two peaks (Equation \ref{eq:energy2a}), we see that the single peak is a lower-energy solution.  This suggests that we might also see a merging of the two peaks, as well as the mass outside the peaks tending to zero.  

Indeed, in our numerical experiments, we saw a merging of peaks except in the special case where $x_0=0.5$, so that the initial peaks are evenly-spaced.  Figure \ref{fig:nle_twin_peaks}a,b shows an example where $x_0=0.5$ but $\epsilon>0$.  Here two peaks remain but the the mass outside those two peaks is absorbed into the peaks over time.   Figure \ref{fig:nle_twin_peaks}c gives an example of peak merging for $x_0<0.5$ whilst Figure \ref{fig:nle_twin_peaks}d shows how the time it takes for peaks to merge increases dramatically as $x_0$ increases towards $x_0=0.5$.  Here, the time to merge is defined as the time at which the centre of the two initial peaks drops below 0.1.  Whilst this is a rather arbitrary definition, other definitions lead to similar trends.

Notice here that the energy analysis does not give direct insight into why merging does not happen for  $x_0=0.5$.  Instead, we turn to physical intuition: the fact that peaks are evenly-spaced means that there is no `preferred' direction for them to move in order to coalesce.  Therefore they remain as two peaks.

\subsection{Peaks of differing heights}
\label{sec:diff_heights}

In Section \ref{sec:same_height}, we examined situations where there are two  peaks with precisely equal height, finding that both peaks persisted indefinitely when they are evenly-spaced.  However, we have already seen in Figure \ref{fig:nle_rand_ic} that when peaks are of different heights, the smaller one can shrink over time, whereas the larger one grows.  If this continues indefinitely, the smaller peak could decay completely and only one peak would remain, although it might take a long time for this to happen. 

Here, we seek to explain this phenomenon using our energy approach, ascertaining whether we should always expect a smaller peak to end up decaying to zero, or whether there are situations where two peaks remain.  To this end, we examine steady state solutions with the following functional form
\begin{align}
	u_*(x)=\begin{cases}
		c_{A}\left[1+\cos\left(\frac{(x+x_0)\sqrt{2}}{\sigma}\right)\right], & \mbox{if }x \in \left(-x_0-\frac{\pi\sigma}{\sqrt{2}},-x_0+\frac{\pi\sigma}{\sqrt{2}}\right) \\
		c_{B}\left[1+\cos\left(\frac{(x-x_0)\sqrt{2}}{\sigma}\right)\right], & \mbox{if }x \in \left(x_0-\frac{\pi\sigma}{\sqrt{2}},x_0+\frac{\pi\sigma}{\sqrt{2}}\right) \\
		0,  & \mbox{otherwise.}
	\end{cases}
	\label{eq:two_humps3}
\end{align}
In this case, we can use Equation (\ref{eq:pt}) to calculate 
\begin{align}
	\label{eq:ca}
	c_{A}=\frac{p-2L\epsilon}{\sqrt{2}\pi\sigma}-c_B.
\end{align}
\begin{figure}[t]
	\begin{center}
		\includegraphics[width=\columnwidth]{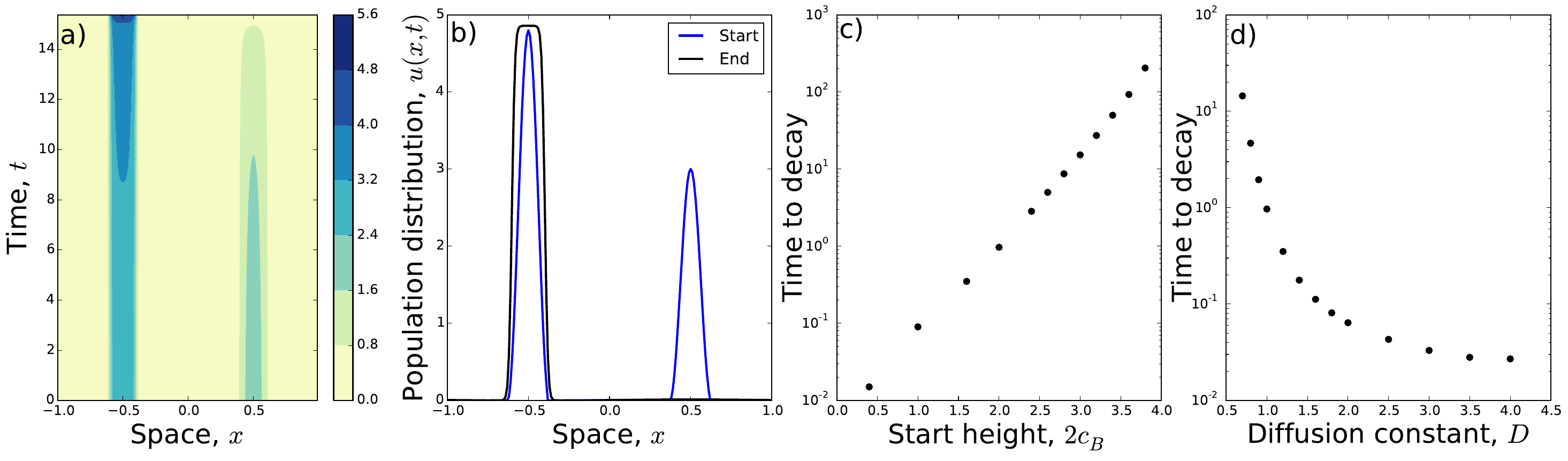}
	\end{center}
	\caption{Panel (a) shows a numerical solution of Equation (\ref{eq:agg2}) with initial condition given by Equation (\ref{eq:two_humps3}) with $c_B=1.5$.		Panel (b) gives snapshots of the initial and final distributions.  Notice that the smaller peak has decayed almost completely by $t \approx 15$.  Panel (c) is constructed from numerical solutions of Equation (\ref{eq:agg2}) with initial condition given by Equation (\ref{eq:two_humps3}) but with $c_B$ taking a variety of values, giving different start heights for the smaller peak (note that the start height is $2c_B$).  Panels (c) and (d) plot the time it takes for the smaller peak to decay to a maximum height of less than 0.1.  This increases exponentially as a function of the start height, explaining the appearance of long-transient multi-peaked solutions to Equation (\ref{eq:agg2}) (Panel c). Conversely, the decay time decreases as $D$ is increased, showing how diffusion can speed up decay of the smaller peak (Panel d). In panels (a-c), $D=1$.  In panel (d), $c_B=1$. In all panels, $\gamma=10$ and $K$ is a top-hat kernel (Equation \ref{eq:Kth}) with $\delta=0.1$. The value of $c_A$ is determined by Equation (\ref{eq:ca}).}
	\label{fig:nle_tdc_plot}
\end{figure}
We see immediately that, in order for $c_A$ and $c_B$ to be non-negative, we must have $c_A,c_B \in \left[0,\frac{p}{\sqrt{2}\pi\sigma}\right]$. A direct calculation using the definition of $E[u]$ from Equation (\ref{eq:energy}) leads to
\begin{align}
	\label{eq:energy2c}
	E[u_*]=-2\sqrt{2}\pi\sigma c_B^2+2pc_B-\frac{p^2}{\pi\sigma\sqrt{2}}.
\end{align}
This is a negative quadratic in $c_B$ with critical point at $c_B=c_A=\frac{p}{2\sqrt{2}\pi\sigma}$.  Therefore the energy minima occur either when $c_B=0$, $c_A=\frac{p}{\sqrt{2}\pi\sigma}$ or  $c_A=0$, $c_B=\frac{p}{\sqrt{2}\pi\sigma}$.  In other words, they occur when there is just one peak.  Consequently, away from the critical point where $c_A=c_B$, we would expect the smaller peak to slowly decay to zero over time, leaving just one peak.  Indeed, this is what we see in numerical solutions of Equation (\ref{eq:agg2}) (e.g. Figure \ref{fig:nle_tdc_plot}a,b).  However, the time it takes for the smaller peak to decay can be very large (Figure \ref{fig:nle_tdc_plot}c).  This is exacerbated by decreasing the diffusion constant, $D$ (Figure \ref{fig:nle_tdc_plot}d).  Here, numerics hint that, as $D \rightarrow 0$, the decay time may tend to infinity, meaning that the second peak may persist indefinitely if there is no diffusion to allow the smaller to seep into the larger.  

\subsection{Including population growth}

So far, we have studied a system where the population size remains constant.  This assumes that there are negligible births or deaths on the timescales that we are studying.  Our focus has been on examining the difference between long transients and asymptotic solutions.  However, in any real biological system, the effect of births and deaths will become non-negligible at some point in time.  Therefore there is a limit to which transient solutions in these systems are biologically realistic: if the transients persist for too long, it will become necessary to account for the effect of births and deaths in any biologically meaningful model.

We therefore examine the extent to which incorporating growth might enable a second peak to persist, by solving the following equation numerically
\begin{align}
	\label{eq:agg_growth}
	\frac{\partial u}{\partial t}=D\frac{\partial^2 u}{\partial x^2}-\gamma\frac{\partial}{\partial x}\left[u\frac{\partial}{\partial x}(K\ast u)\right]+ru\left(1-\frac{u}{K}\right),
\end{align}
with initial conditions given by Equation (\ref{eq:two_humps3}). 

\begin{figure}[h!]
	\begin{center}
		\includegraphics[width=\columnwidth]{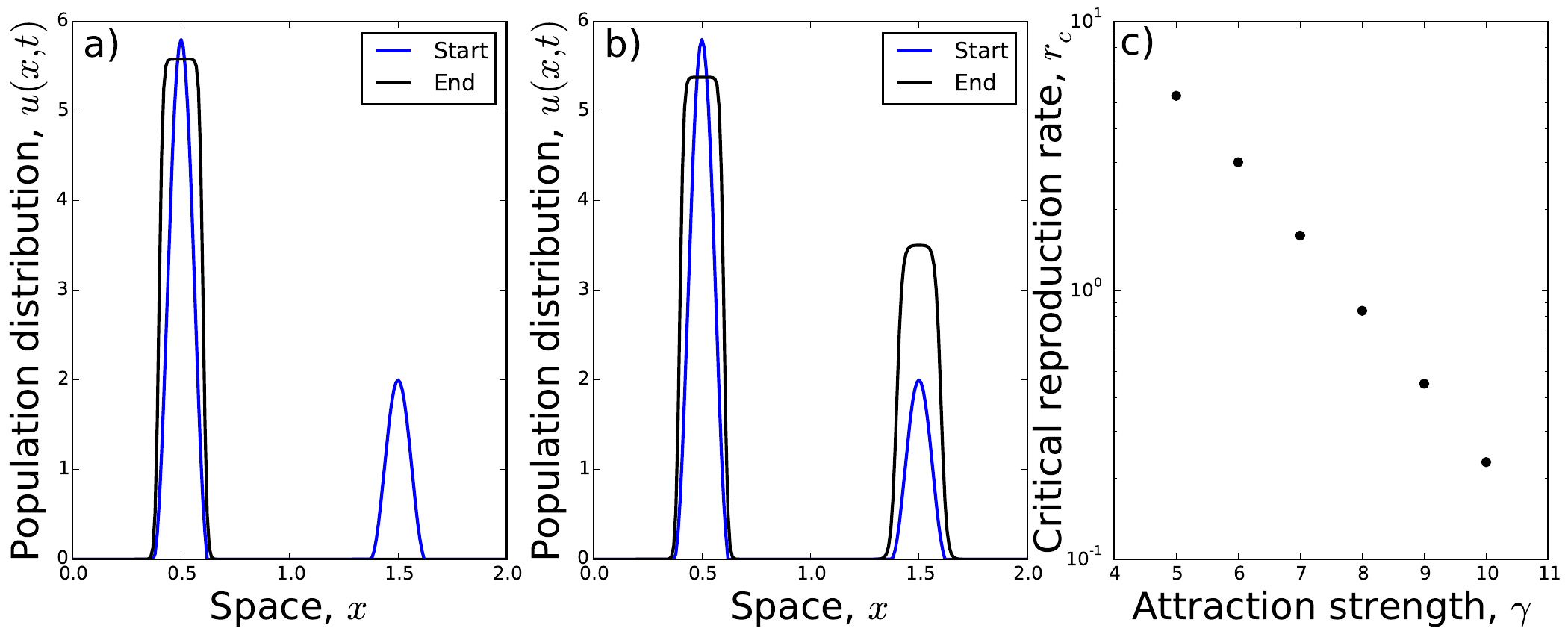}
	\end{center}
	\caption{{\bf Effect of growth parameter.}  Panels (a) and (b) show the initial condition (blue) and solution at time $t=10$ (black) where the parameters are $\gamma=10$, $D=1$, and $K=5$.  In Panel (a), $r=0.23$ whereas Panel (b) has $r=0.24$.  This demonstrates a transition in long-term patterns, whereby the smaller peak decays for $r\leq 0.23$ but grows for $r\geq 0.24$.  Panel (c) shows how this transition point, $r_c$, decreases exponentially as the strength of attraction, $\gamma$, increases.}
	\label{fig:nle_growth_plot}
\end{figure}
Depending upon the values of $\gamma$, $D$, $K$, and $c_B$, we found that there is a critical value $r=r_c$ above which the second hump persists, and below which it decays.  Figure \ref{fig:nle_growth_plot}a,b shows this in the case $\gamma=10$, $D=1$, $K=5$, $c_B=1$, whereby $r_c\approx 0.23$.  Figure \ref{fig:nle_growth_plot}c demonstrates how $r_c$ depends upon the aggregation strength $\gamma$: the greater the aggregation strength, the higher the required growth rate to enable a second peak to persist.

\section{Discussion}

Distinguishing between asymptotic solutions and long transients in numerical PDEs is a thorny issue, with perhaps no one-size-fits-all solution.  Typically, researchers decide that a solution has reached an asymptotically-stable state when some measure (e.g. the change in $L^p$ norm for some $p \in [1,\infty]$) is below a small threshold value (see e.g.  \citet{burger2014stationary,giunta2022local,schlichting2022scharfetter}).  However, this means that if transient solutions are changing slower than this threshold value then they will be mistaken for asymptotically-stable solutions.  Therefore it is valuable to have some analytic insight to guide the user as to whether the solution is (or is likely to be) a long transient or an asymptotically-stable solution.

Here, we have provided such a deductive technique for the aggregation-diffusion equation in Equation (\ref{eq:agg}).  Rather than studying this equation directly, we instead study an approximation given in Equation (\ref{eq:agg2}).  This approximate formulation is simple enough to solve for steady state solutions.  It also possesses an energy functional, which allows us to search for local minimum energy solutions amongst the steady state solutions, an approach employed successfully in a previous multi-species study \citep{giunta2022detecting}.  Our hypotheses are first that these local minimum energy solutions are stable solutions to Equation (\ref{eq:agg2}), whereas other steady states are not; and second that this categorisation carries over to the steady states of Equation (\ref{eq:agg}).  In the examples we tested, numerical experiments confirmed these hypotheses, with the sole exception of twin-peaked solutions where the peaks are of identical height and evenly-spaced.  We therefore conclude that this method is a useful way for guiding users (i.e. those wanting to solving Equation \ref{eq:agg} numerically) as to whether a solution they are observing is likely to be stable or not, whilst also recommending that they verify these calculations up with numerical experiments.

Regarding the examples we tested, we found two main results: first, that stable aggregations are likely to resemble compactly-supported solutions, rather than being non-zero everywhere; second, that multi-peaked solutions will always be transient unless either $D=0$ or the peaks are precisely the same height and evenly-spaced.  In addition to these central messages, further numerical investigations revealed that these twin-peaked transient solutions can be arbitrarily long-lived if the peaks are arbitrarily close to being evenly-spaced (Figure \ref{fig:nle_twin_peaks}) and the heights of these peaks are arbitrarily similar (Figure \ref{fig:nle_tdc_plot}).  

That said, the consideration of very long transients in a model that operates on timescales where births and deaths are negligible is not terribly realistic, so we also examined the effect of adding a small amount of (logistic) growth.  We found that arbitrarily small amounts of growth will not stop the smaller peak from decaying.  However, there appears to be a critical growth rate, dependent upon the model parameters, below which the smaller peak will decay and above which it will grow (Figure \ref{fig:nle_growth_plot}).  Therefore, if long transients appear when using Equation (\ref{eq:agg}) to model biological aggregation, it is valuable to think about the effect of net reproductive rate in the system being modelled, and whether this is sufficient to arrest the decay of the smaller peak.

Whilst our principal equation of interest is Equation (\ref{eq:agg}), it is worth noting that our approximate analytic techniques can also be applied to various other Equations.  For example, the cell adhesion equations introduced in \citet{armstrong2006continuum},  have a very similar functional form that can usually be formally related to  Equation (\ref{eq:agg}) or  modifications thereof \citep{painter2023biological}.  
Chemotaxis equations are also somewhat similar to Equation (\ref{eq:agg}), but here the non-local self-interaction is replaced with a diffusing chemical.  The organisms interact with the chemical rather than directly with one another.  It turns out that the resulting models are equivalent to a type of aggregation-diffusion equation with advection that is nonlocal in both space and time \citep{shi2021spatial}. This contrasts with Equation (\ref{eq:agg}), which is nonlocal in space alone.  However, similar patterns are observed in these systems, including long-transient multi-peaked solutions similar to those studied here \citep{potapov2005metastability}.  We also note that the moment closure we use in Equation (\ref{eq:agg}) leads to a fourth-order PDE quite similar in nature to the Cahn-Hilliard equation \citep{novick2008cahn}, for which there is a long history of studies on metastability \citep{bates1994metastable, reyna1995metastable, scholtes2018metastability}.  

Finally, it is worth noting that the particular version of the aggregation-diffusion equation that we study involves linear diffusion.  However, there is also interest in the nonlinear case, particularly where the diffusion is quadratic, replacing $u_{xx}$ with $(u^2)_{xx}=2(uu_x)x$ in Equation (\ref{eq:agg}).  An advantage of this formulation is that Equation (\ref{eq:agg}) has the form $u_t=[u(D-\gamma K\ast u)_x]_x$, making it amenable to a analysis without taking the limit $D/\gamma \rightarrow 0$.  This fact has been exploited, for example, by \citet{ellefsen2021equilibrium,carrillo2018zoology}.   However, here we have chosen to focus on linear diffusion is important to study as it often arises naturally from models of organism movement \citep{armstrong2006continuum, potts2020parametrizing, painter2023biological}.  Future work on the nonlinear case could reveal analytic insights about the effect of $D$ vs. $\gamma$ on asymptotic patterns, which we were only able to examine numerically in this study. 

\section*{Acknowledgments}

JRP acknowledges support of Engineering  and  Physical  Sciences  Research  Council (EPSRC) grant EP/V002988/1. KJP is a member of INdAM-GNFM and acknowledges departmental funding through the `MIUR-Dipartimento di Eccellenza' programme.





\bibliography{nle_refs}

\end{document}